\documentclass{article}
\usepackage[utf8]{inputenc}

\usepackage{hyperref,mathrsfs,enumitem,amsmath, amssymb, amsthm, graphicx,color}

\DeclareMathOperator{\Log}{Log}

\newcommand{\vertiii}[1]{{\left\vert\kern-0.25ex\left\vert\kern-0.25ex\le
ft\vert #1 
    \right\vert\kern-0.25ex\right\vert\kern-0.25ex\right\vert}}

\newcommand{\C}{\mathbb{C}}

\newcommand{\D}{\mathbb{D}}

\newcommand{\ds}{\displaystyle}
\usepackage{mathtools}

\DeclarePairedDelimiter\floor{\lfloor}{\rfloor}
\newtheorem{theorem}{Theorem}
\newtheorem{lemma}[theorem]{Lemma}

\newtheorem{proposition}[theorem]{Proposition}

\begin{document}

\title{Linear Factorization of Hypercyclic Functions for Differential Operators}


%\fntext[myfootnote]{Present address: Department of Mathematics, Statistics, and Computer Science, St. Olaf College, Northfield, MN 55057 }

\author{Kit C. Chan\footnote{Department of Mathematics and Statistics, Bowling Green State University, Bowling Green, OH 43403}, Jakob Hofstad\footnote{Department of Mathematics, Statistics and Computer Science, St. Olaf College, Northfield, MN 55057}, David Walmsley\textsuperscript{$\dagger$}}

\date{}
\maketitle

%\address[mymainaddress]{Department of Mathematics and Statistics, Bowling Green State University, Bowling Green, OH 43403}
%\address[mysecondaryaddress]{Department of Mathematics, Statistics and Computer Science, St. Olaf College, Northfield, MN 55057}

\begin{abstract}
    On the Fr\'{e}chet space of entire functions $H(\C)$, we show that every nonscalar continuous linear operator $L:H(\C)\to H(\C)$ which commutes with differentiation has a hypercyclic vector $f(z)$ in the form of the infinite product of linear polynomials:  
\[
f(z) = \prod_{j=1}^\infty \, \left( 1-\frac{z}{a_j}\right),
\]
where each $a_j$ is a nonzero complex number.
\end{abstract}


\section{Introduction}\label{Introduction}
A sequence of continuous maps $\{T_n: X \to X|n \geq 1\}$ on a separable metrizable topological space X is said to be {\it universal} if there is an 
element $x$ in $X$ such that the set $\{T_nx : n \geq 1\}$ is dense in $X$. Such an element $x$ is called a {\it universal element} for the sequence 
$\{T_n\}$. In the special case when there is a map $T$ such that each member $T_n$ of the universal sequence is given by $T_n = T^n$, then we say that the map $T$ is {\it hypercyclic}, and a universal element for $\{T^n\}$ is called a {\it hypercyclic element } for $T$.

Some of the early examples of universality were exhibited on the vector space $H(G)$ of all analytic functions on a region $G$ of the complex plane 
$\mathbb{C}$. Given the compact-open topology, $H(G)$ becomes a Fr\'echet space, in which a sequence $\{f_n\}$ converges to $f$ if and only if 
$f_n \to f$ uniformly on compact subsets of $G$. When the region $G = \mathbb{C}$, Birkhoff  \cite{Birkhoff}  showed in 1929 that the translation operator  
$T:f(z) \mapsto f(z+1)$ is hypercyclic on $H(\C)$. In 1952, MacLane \cite{MacLane} showed that the differentiation operator $D: f(z) \mapsto f'(z)$ is hypercyclic on $H(\C)$. 
These two are classical examples of hypercyclicity on the space  $H(\C)$ of entire functions.

Extending the results of Birkhoff and MacLane, Godefroy and Shapiro \cite[Theorem 5.1]{GodefroyShapiro}  proved in 1991 that every continuous 
nonscalar linear operator $L:H(\mathbb{C}) \to H(\mathbb{C})$ commuting with the differentiation operator $D$ is hypercyclic.  Furthermore, they \cite[Theorem 5.2]{GodefroyShapiro} 
showed that every continuous linear operator $L:H(\mathbb{C}) \to H(\mathbb{C})$ commuting with $D$  is of the form $T = \varphi (D)$ where 
$\varphi$ is an entire function of exponential type. Such operators $\varphi (D)$ are called {\it differential operators}.  For basic properties of differential operators, we refer the reader 
to the book of Grosse-Erdmann and Peris Manguillot \cite[Section 4.2]{LinChaos}.  

However, one question remains unanswered. That is, can we exhibit a particular 
hypercyclic vector of any differential operator $\varphi(D)$, with certain structure?  To answer this question, we show in Theorem \ref{HypThm} below that every such operator $\varphi (D)$ which is not a scalar multiple of the identity has 
a hypercyclic vector $f(z)$ that is in the form of the infinite product of linear polynomials:  
\begin{eqnarray}\label{goal}
f(z) = \prod_{j=1}^\infty \, \bigg( 1-\frac{z}{a_j}\bigg),
\end{eqnarray}
where each $a_j$ is a nonzero complex number.

The rate at which the zeros of an entire function $f$ tend to $\infty$ is closely related to the
exponential growth of $f$ and to the structure of the factors in the Weierstrass product expansion of $f$.  The growth of $\varphi(D)$-hypercyclic functions has been studied by several authors. We refer the reader to the paper by Bernal-Gonz\'{a}lez and Bonilla \cite{BernalBonilla}, and the references therein. Functions of the above form in (\ref{goal}) were used by Chan and Shapiro \cite{ChanShapiro} in showing that the translation operator on a Hilbert space of entire functions
of slow growth is hypercyclic, however no specific form of a hypercyclic vector was given.  Their techniques cannot be easily adapted to work for the 
differentiation operator $D$, because repeated differentiation of  an infinite  product results in a very complicated formula.   This differentiation problem is made much worse, if we do that for a more general operator 
$\varphi(D)$ that we consider in the present paper.

The main result Theorem \ref{HypThm} of the present paper is analogous to a well-known result for the case of the open unit disk $\mathbb{D}$.  To explain that, let $\psi_n: \mathbb{D} \to \mathbb{D}$ be 
a sequence of analytic self-maps of $\mathbb{D}$, and consider the sequence of composition operators $T_n: H(\mathbb{D} ) \to H(\mathbb{D})$ defined by $T_nf = f\circ \psi_n.$   In 1941, Seidel and Walsh \cite{SeidelWalsh} gave the first 
example of a universal sequence of composition operators $\{T_n\}$ on  $H(\mathbb{D})$. Heins, a student of Walsh's, focused on the multiplicative semisubgroup $X = $ Ball\,$H^\infty(\mathbb{D})$  $=$   $ \{f \in H(\mathbb{D}):$   $|f(z)| \leq 1,$ for all  $z\in \mathbb{D}\}$ of $H(\mathbb{D})$ in \cite{Heins}.  He 
 showed in 1954 that for any universal sequence $\{T_n\}$ on $H(\D)$ given in \cite{SeidelWalsh},  there is a Blaschke product $b$ that is universal 
 for $\{T_n\}$ on $X$.   Observe that a  Blaschke product is a convergent product of automorphisms of $\mathbb{D}$. Since  linear polynomials are automorphisms 
 of $\mathbb{C}$,  the
    infinite products in $H(\C)$ given in (\ref{goal}) above are analogous to the Blaschke products in $H(\D)$.  

The result of Heins was generalized by Bayart, Gorkin, Grivaux and Mortini \cite{BGGM} to the following result:  For a sequence of analytic self-maps $\psi_n: \D \to \D$ such that $\limsup_{n\to \infty} |\psi_n(0)|=1$, the corresponding
sequence of composition operators $\{T_n\}$ is universal on $X$ if and only if there is a Blaschke product $b$ that is universal for $\{T_n\}$.  
Their result was recently given a simpler proof  by Chan and Walmsley \cite{ChanWalmsley} using a constructive approach to generalize  the  Universal Criterion to a semigroup setting.  This constructive approach allowed them to exhibit a Blaschke product as a universal element. 
The original form of  Universal Criterion given by Gethner and Shapiro \cite{GethnerShapiro} can be applied to prove the aforementioned results of 
Birkhoff   \cite{Birkhoff}, MacLane \cite{MacLane}, and Godefroy and Shapiro \cite[Theorem 5.1]{GodefroyShapiro}, but it does not provide a specific form of a universal element.  

Returning to the results of the present paper, we remark that due to the complex structure of $\varphi(D)$, where $\varphi$ is a nonconstant entire function of exponential type and $D$ is the differentiation operator,  the constructive approach of Chan and Walmsley \cite{ChanWalmsley}  does not work immediately to provide a specific form of 
hypercyclic vector for $\varphi(D)$.   In fact, we cannot use 
the techniques of employing the Universality Criterion in any way.  We basically have to construct a hypercyclic vector of the form given by (\ref{goal}).  
Thus our techniques are very different from those commonly found in the area of hypercyclicity and universality. 
 
In Section \ref{Estimates on Polynomials} below, we obtain four technical lemmas for  making estimations on polynomials.   Then in Section \ref{QuotientPolys},  we focus on quotient polynomials, which play an important 
role in showing the hypercyclicity of the function of the form given by (\ref{goal}).
Lastly, in Section \ref{MainResult}, we prove our main result Theorem \ref{HypThm} using the technical results in Sections \ref{Estimates on Polynomials} and \ref{QuotientPolys}.

%%%%%%%%%%%%%%%%%%%%%%%%%%%%%%%%%%%
%%%%% ESTIMATES ON POLYNOMIALS %%%%
%%%%%%%%%%%%%%%%%%%%%%%%%%%%%%%%%%%
\section{Estimates on Polynomials}\label{Estimates on Polynomials}

In this section we prove a few lemmas on polynomials that facilitate our estimations in the proof of the main result Theorem \ref{HypThm}. Let $D$ be the differentiation operator and $I$ be the identity operator on $H(\C)$, and let $\mathcal{P}$ be the collection of complex polynomials in $H(\C)$. If $R>0$ and $f\in H(\C)$, let $\|f\|_R = \max_{|z|\leq R} |f(z)|$.  We first provide a description of a right inverse on $\mathcal{P}$ for certain types of differential operators.

%%%%%%%%%%%%%%%%%%%%%%%%%
%%%%% INVERSE LEMMA %%%%%
%%%%%%%%%%%%%%%%%%%%%%%%%

\begin{lemma}\label{InverseLemma}   Let $p \in \mathcal{P}$ and let $m$ be the degree of $p$. If  $T=\varphi(D)$ is a differential operator for some entire function $\varphi$ of exponential type with $\varphi(0)\not =0$, then there exists a mapping  $S: \mathcal{P} \to \mathcal{P}$ such that $Sp$ has degree $m$ and $T S p = p$. Moreover, if $c_{n,i}$ denotes the coefficient of $z^i$ of $S^n p$, then there exists a constant $C>1$ that does not depend on $n$ such that
\[|c_{n,i}| < C^n. \]
\end{lemma} 

\begin{proof}
We first write $ T = \displaystyle\sum\limits_{j=0}^{\infty} a_j D^j$, where $a_0\not = 0$.  Since $D^n p = 0$ for any integer $n>m$, we have $T p =(a_0 I + a_1 D + \cdots + a_m D^m)p$.  If $\alpha_i \in \C$ are the zeros of the polynomial $a_0+a_1 z + \cdots + a_m z^m$, repeated according to multiplicity, then since $a_0\not =0$, none of the $ \alpha_i $ equal zero and we can write 
\begin{align}\notag
Tp = a_0(I-D/\alpha_1)(I-D/\alpha_2)(I-D/\alpha_3)\cdots(I-D/\alpha_{m-1})(I-D/\alpha_m)p
\end{align}

To find a formula for the mapping $S$ in the statement of our lemma, we first find a right inverse $S_i$ for each factor $I-D/\alpha_i$. For that we observe
\begin{align*}
p &= (I - D^{m+1}/\alpha_i^{m+1})p \\
    &= (I - D /\alpha_i)(I+(D/\alpha_i)+(D/\alpha_i)^2+ \dots + (D/\alpha_i)^m)p.
\end{align*}
Thus define $S_i p = (I+(D/\alpha_i)+(D/\alpha_i)^2+ \dots + (D/\alpha_i)^m)p$, which is a polynomial of degree $m$. We then define $Sp$ as $Sp=\frac{1}{a_0}S_1\cdots S_m p$.  Thus $Sp$ has degree $m$ and $TSp=p$.

By writing the formula for $Sp$ as
\begin{align}\notag
S p = \frac{1}{a_0}\left( \displaystyle\prod\limits_{i=1}^{m} (I + D/\alpha_i + D^2/\alpha_i^2+\dots+D^m/ \alpha_i^m) \right) p,
\end{align}
we now proceed to obtain bounds on the coefficients of the polynomial $S^n p$. Since $p$ has degree $m$, we first write
\begin{align}\label{Snpeqn}
S^n p= \frac{1}{a_0^n}(b_{0,n}I + b_{1,n}D+b_{2,n}D^2+\dots+b_{m,n}D^m) p.
\end{align}
Let $r=\max\{1,|\alpha_1|^{-1},|\alpha_2|^{-1},\ldots,|\alpha_m|^{-1}\}$ and let $C(mn,i)$ be the coefficient of $y^i$ in the expansion of $(1+y+y^2+y^3+...)^{mn}$. By multiplying out $S^n$, we have that $|b_{i,n}| \leq r^i C(mn,i)$.  Since $(1+y+y^2+\dots)^{mn} = 1/(1-y)^{mn}$ for $y \in (-1,1)$, and the Taylor Series for $1/(1-y)^{mn}$ is 
\[1+mny+\frac{mn(mn+1)y^2}{2!}+\frac{mn(mn+1)(mn+2)y^2}{3!}+\dots,\]
we have that 
\[C(mn,i) =\binom{mn+i-1}{i} =  \left( \frac{mn}{1}\right) \left( \frac{mn+1}{2}\right)\cdots \left( \frac{mn+i-1}{i}\right) \leq (mn)^i,\]
which implies that 
\begin{align}\label{b_i,n}
    |b_{i,n}|\leq r^i C(mn,i) \leq (rmn)^i.
\end{align}

\sloppy With this in mind, we now provide a bound for the coefficient $c_{n,i}$ of $z^i$ in the polynomial $S^n p$. Let $p(z) = p_m z^m + p_{m-1} z^{m-1}+\dots+p_0 $ and $\mu=\max\{1,|p_0|,|p_1|,\ldots,|p_m|\}$.  Then by (\ref{Snpeqn}), the coefficient $c_{n,i}$ satisfies
\begin{align}\notag
    |c_{n,i}| &= \left| \frac{1}{a_0^n}\displaystyle\sum\limits_{j=i}^{m} b_{j-i,n} \ (D^{j-i}p_j z^j)|_{z=1} \right|\\ \notag
    & = \left| \frac{1}{a_0^n} \displaystyle\sum\limits_{j=i}^{m} b_{j-i,n} \ p_j \frac{j!}{i!} \right|\\ \notag
    & \leq \left| \frac{1}{a_0} \right|^n \displaystyle\sum\limits_{j=i}^{m} \left| (rmn)^{j-i} p_j \frac{j!}{i!} \right| \text{ by (\ref{b_i,n})}\\
    & < \left| \frac{1}{a_0} \right|^n (m+1)(rm)^m \mu m! n^m. \label{coefficient}
\end{align}
Now choose a constant $\gamma$ such that $\gamma>\max\{1,|a_0|^{-1},(m+1)!(rm)^m \mu, e^m\}$.  Then $\ds \frac{\ln \gamma}{m}>1$, which implies that $\ln n < n < n \ds \frac{\ln \gamma}{m}$.  Thus $\ln(n^m)<\ln(\gamma^n)$, which yields that $\ds \frac{n^m}{\gamma^n}<1$.  Hence
\begin{align}
    \frac{|a_0|^{-n}}{\gamma^n} \frac{(m+1)!(rm)^m \mu}{\gamma} \frac{n^m}{\gamma^n} <1. \label{constant}
\end{align}
Now if we let $C = \gamma^3$, then one can easily verify that $C$ does not depend on $n$ and that $C>1$. An application of inequality (\ref{constant}) to inequality (\ref{coefficient}) yields that
\[ |c_{n,i}| < \gamma^{2n+1} \leq \gamma^{3n}=C^n.\]
\end{proof} 

For a polynomial $p$ and differential operator $T$, our next lemma is to estimate the growth of $T^n p$ on a disk $\{z: |z|\leq R\}$.  It is crucial for our construction of a hypercyclic function later on.
%%%%%%%%%%%%%%%%%%%%%%%%%%%%%
%%%%% POLY GROWTH LEMMA %%%%%
%%%%%%%%%%%%%%%%%%%%%%%%%%%%%

\begin{lemma}\label{PolyGrowthLemma}
Suppose $T$ is a differential operator. Let $p(z)=p_mz^m+\cdots+p_0$ be a polynomial of degree $m$, let $\mu=\max\{|p_j|: 0\leq j\leq m\}$, and let $R>0$ be given.  There exist constants $\alpha >1$ and $\beta >1$ that depend only on $T$ such that $\|T^n p\|_R < \mu (m+1) \alpha^n (n\beta +R)^m$.
\end{lemma}

\begin{proof}
Let $T=\varphi(D)$, where $\varphi(z)=\sum_{j=0}^\infty a_j z^j$ is an entire function of exponential type. By \cite[Lemma 4.18]{LinChaos}, there exist constants $\alpha>1$ and $\beta>1$ such that 
\begin{align*}
    |a_j|\leq \frac{\alpha \beta^j}{j!}, (j=0,1,2,\cdots).
\end{align*}
Let $(\varphi(z))^n = \sum_{j=0}^\infty \gamma_{n,j} z^j$. We first show by induction that 
\begin{align}
    |\gamma_{n,j}| \leq \frac{\alpha^n (n\beta)^j}{j!}.\label{betaEqn}
\end{align}

The case $n=1$ is trivial.  Now suppose (\ref{betaEqn}) holds for $1,\ldots, n$.  Then
\begin{align*}
    (\varphi(z))^{n+1} = (\varphi(z))^n(\varphi(z)) &= \left( \sum\limits_{j=0}^{\infty} \gamma_{n,j} z^j \right) \left( \sum\limits_{k=0}^{\infty} a_k z^k \right) \\
    & = \sum\limits_{j=0}^{\infty}\left( \sum\limits_{k=0}^{j} \gamma_{n,j-k} \, a_{k} \right) z^j.
\end{align*}
Hence
\begin{align*}
    |\gamma_{n+1,j}| &= \left| \sum\limits_{k=0}^{j} \gamma_{n,j-k} \, a_{k} \right| \\
    & \leq\sum\limits_{k=0}^{j} |\gamma_{n,j-k}| \, |a_{k}| \\
    & <\sum\limits_{k=0}^{j} \left( \frac{\alpha^n (n\beta)^{j-k}}{(j-k)!} \right) \, \left( \frac{\alpha\beta^k}{k!} \right)  \text{ by the Induction Hypothesis}\\
    & = \left( \frac{\alpha^{n+1} \beta^j}{j!} \right) \sum\limits_{k=0}^{j}  \,  n^{j-k} \binom{j}{k} \\&=
\frac{\alpha^{n+1} (\beta (n+1))^j}{j!}  \text{ by the Binomial Theorem.}
\end{align*}
This verifies (\ref{betaEqn}) and the induction is complete.

Since $T^n=(\varphi(D))^n$, we have for any monomial $z^k$ that 
\begin{align*}
    \|T^n(z^k)\|_R & = \left\|\sum_{j=0}^k \gamma_{n,j}D^j (z^k) \right\|_R\\
    & \leq \sum_{j=0}^k |\gamma_{n,j}| R^{k-j}\frac{k!}{(k-j)!}\\
    & \leq \sum_{j=0}^k \frac{\alpha^n (n\beta)^j}{j!}  R^{k-j} \frac{k!}{(k-j)!}  \text{ by (\ref{betaEqn}) }\\
    & = \alpha^n \sum_{j=0}^k (n\beta)^j R^{k-j} \binom{k}{j}\\
    & = \alpha^n (n\beta +R)^k.
\end{align*}
Therefore $\|T^n p\|_R \leq \sum_{k=0}^m \left\| p_k T^n(z^k)\right\|_R  \leq \mu (m+1) \alpha^n (n\beta +R)^m$.
\end{proof}

%%%%%%%%%%%%%%%%%%%%
%%%%% r Lemma %%%%%%
%%%%%%%%%%%%%%%%%%%%

Using Lemmas \ref{InverseLemma} and \ref{PolyGrowthLemma}, we now construct a sequence of functions with a desirable limit property that facilitates the proof of our main theorem.

\begin{lemma}\label{rLemma}
Suppose $T = \varphi (D)$ is a differential operator for some non-constant entire function $\varphi$ of exponential type with $\varphi (0) \not= 0$.  Let $S:\mathcal{P}\to \mathcal{P}$ be the right inverse given by Lemma \ref{InverseLemma}. There exists a nonzero $r\in \C$ with the following property: if $f$ is a polynomial of degree $m$, $p$ is a polynomial of degree less than $m$, and $n> m-1$, then as $n \to \infty$,

\[T^n \left( (f(z) - S^n p(z))  \sum_{i=0}^{n-m-1} \frac{(-rz)^i}{i!}  \sum_{j=0}^{\floor*{ n^{1.2}}}  \frac{(rz)^j}{j!} \right) \to 0\]
uniformly on compact subsets of $\C$.
\end{lemma}

\begin{proof}
We first show there exists a nonzero $r\in \C$ such that
\begin{align}\label{e^r}
    T^n \left( (f(z) - S^n p(z)) \left( \sum_{i=0}^{n-m-1} \frac{(-rz)^i}{i!} \right)e^{rz}\right) \to 0.
\end{align}
Since $\varphi(0)\not =0$, our proof of (\ref{e^r}) is divided into two cases depending on whether $\varphi(z)$ has a zero at $w\not=0$ or has no zero. For the first case that $\varphi(w)=0$ for some $w\not =0$, let $r=w$ and write $\varphi(z)=\psi(z)(1-z/r)$ where $\psi$ is an entire function, so that $T=\psi(D)(I-D/r)$. If $h$ is a polynomial, one can quickly check that $(I - D/r)(e^{rz} h(z)) = (-1/r)e^{rz} D(h(z))$.
Furthermore, if $h$ is a polynomial of degree less than $n$, then by an inductive argument we obtain $(I-D/r)^n (h(z)e^{rz})=0$, which yields that
\begin{align}\label{herz}
    T^n (h(z)e^{rz}) = 0.
\end{align}
Since by Lemma \ref{InverseLemma} the degree of $S^n p$ equals the degree of $p$, the function $(f(z) - S^n p(z)) \left( \sum_{i=0}^{n-m-1} \frac{(-rz)^i}{i!} \right)$ is a polynomial of degree $n-1$. Thus, in this case, by using (\ref{herz}) we see that $T^n \left( (f(z) - S^n p(z)) \left( \sum_{i=0}^{n-m-1} \frac{(-rz)^i}{i!} \right)e^{rz}\right) = 0$, so the limit in (\ref{e^r}) holds.

For the second case that $\varphi(z)$ has no zero, we write $\varphi (z) = e^{g(z)}$ for some entire function $g(z)$.  Since $\varphi (z)$ is of exponential type, there must exist $\gamma, \sigma \geq 0$ such that $|\varphi(z)| < \gamma e^{\sigma|z|}$ for all $z \in \C$. Thus by \cite[Thm 4.14.3]{Holland}, we have $g(z)=az+b$ for some $a,b\in \mathbb{C}$.  Since $T$ is a nonscalar operator, $a\not =0$. Thus if we write $\lambda = e^b$ then $T$ has the form $T=\lambda e^{aD}$.

Using the Taylor series expansion at $z$, it is easily seen that $e^{aD}=T_a$, where $T_a: H(\C)\to H(\C)$ is the translation operator given by $T_a h(z) = h(z+a)$.  By Lemma \ref{InverseLemma}, there exists a constant $C> 1$ such that
\begin{align}\label{coefficients}
    \text{the coefficients of the polynomial $f-S^np$ are bounded above by $C^n$.}
\end{align}
Hence by Lemma \ref{PolyGrowthLemma}, for any $R>0$, there are positive constants $\alpha$ and $\beta$ such that
\begin{align} \label{alphabeta}
    \| T_a ^n (f-S^n p) \|_R & \leq C^n (m+1) \alpha^n (n\beta +R)^m\\
    & \leq \kappa^n  \text{ for some  $\kappa>0$. }\label{Ta f-Snp est}
\end{align}

Choose a nonzero $r\in \mathbb{C}$ such that $|r|>1$ and
\begin{align}
     \left| e^{-ra} \right| > \max\{|\lambda rae|\kappa, e^3\} .\label{rEqn}
\end{align}
To continue our argument using the expression $T_a=\sum_{j=0}^\infty \frac{(aD)^j}{j!}$, we let 
\[s_k (z) = T_a^n \left(\frac{(-rz)^k}{k!}\right)= \frac{(-r(z+na))^k}{k!}. \]
Since $|ra|>3$ by (\ref{rEqn}), we see that if $n>|r|R$ and $n>k$, then for  $|z|\leq R$,
\[     \left|\frac{s_k(z)}{s_{k-1}(z)} \right| = \frac{|-r(z+na)|}{k}  >  \frac{|-r(z+na)|}{n} \geq \left|  |ra| -\frac{|rz|}{n} \right| > 2. \]
Thus for $n$, $k$, and $z$ as above,
\begin{align}
    \left|\frac{s_k(z)}{s_n(z)}\right| =\left|\frac{s_k(z)}{s_{k+1}(z)} \frac{s_{k+1}(z)}{s_{k+2}(z)}\cdots \frac{s_{n-1}(z)}{s_n(z)} \right|< \frac{1}{2^{n-k}}. \label{skEqn}
\end{align}
Hence when $n>|r|R$,
\begin{align}\notag
    \left\| T_a ^n\left( \displaystyle\sum\limits_{k=0}^{n-m-1} \frac{(-rz)^k}{k!} \right) \right\|_R & \leq  \sum\limits_{k=0}^{n-m-1} \| s_k\|_R \\ \notag
    & < \sum\limits_{k=0}^{n-m-1} \frac{\|s_n\|_R}{2^{n-k}}  \text{ by (\ref{skEqn}) }\\ \notag
    & < \|s_n\|_R\\ \notag
    & = \left\| \frac{(-r(z+na))^n}{n!} \right\|_R \\ \notag
    & < \left\| \frac{(-r(z+na))^n e^{n}}{n^n} \right\|_R  \text{ because $\frac{n^n}{n!} < e^n$ }\\ \notag
    & \leq \left| (rae)^n \right| \cdot \left\| \left(\frac{z}{na}+1\right)^n  \right\|_R \\ 
    & < \left| (rae)^n \right| e^{R/|a|}. \label{sum bound}
\end{align}

Since $T=\lambda e^{aD}=\lambda T_a$, we have $T(gh)=\lambda T_a(g) T_a(h)$ for any entire functions $g$ and $h$. We compute that
\begin{align}\notag
    & \left\| T^n \left( (f(z) - S^n p(z)) \left( \sum_{i=0}^{n-m-1} \frac{(-rz)^i}{i!} \right)e^{rz}\right) \right\|_R \\ \notag
    & \leq |\lambda|^n \left\| T_a^n (f-S^n p)\right\|_R \left\|T_a ^n\left( \displaystyle\sum\limits_{k=0}^{n-m-1} \frac{(-rz)^k}{k!} \right) \right\|_R \left\| T_a^n (e^{rz}) \right\|_R\\ \notag
    & \leq |\lambda|^n  \kappa^n \left| (rae)^n \right| e^{R/|a|} \|e^{rz+rna}\|_R  \text{ by (\ref{Ta f-Snp est}) and (\ref{sum bound}) }\\ \notag
    & \leq \frac{|\lambda|^n \kappa^n \left| (rae)^n\right| e^{R/|a|} e^{|r|R}}{\left|e^{-rna}\right|} \\
    & =  e^{R/|a|} e^{|r|R} \left|\frac{\lambda rae \kappa}{e^{-ra}}\right|^n. \label{bigTneqn}
\end{align}
Therefore by inequality (\ref{rEqn}), the above expression converges to $0$ as $n\to \infty$, which verifies the limit in (\ref{e^r}).

We are now ready to prove the limit in the conclusion of our lemma. For that, we observe
\[\sum_{j=1}^{\floor*{ n^{1.2}}}  \frac{(rz)^j}{j!}= e^{rz} - \sum_{j=\lceil n^{1.2} \rceil}^\infty \frac{(rz)^j}{j!},\]
and hence by (\ref{bigTneqn}) it suffices to show that 
\begin{align}\label{rlimit}
    T^n \left( (f(z) - S^n p(z))  \sum_{i=0}^{n-m-1} \frac{(-rz)^i}{i!} \sum_{j=\floor*{ n^{1.2}} + 1}^\infty \frac{(rz)^j}{j!} \right) \to 0.
\end{align}
To simplify matters, let $g_n = \ds (f(z) - S^n p(z)) \left( \sum_{i=0}^{n-m-1} \frac{(-rz)^i}{i!} \right)$, which is a polynomial of degree $n-1$.   Let $B=\sup\{|r|^j/j!: 0\leq j <\infty\} $, which is finite. By multiplying out $g_n$ and using (\ref{coefficients}), the absolute values of the coefficients of $g_n$ are bounded above by $(m+1)C^n B$.  Then for any positive integer $j$, the degree of the polynomial $z^j g_n(z)$ is $j+n-1$, and its coefficients are bounded above by $(m+1)C^n B$.  Since $T^n = \lambda^n T_a^n$, the inequality in (\ref{alphabeta}) derived from Lemma 2 implies that
\begin{align}\label{Tnzj}
    \|T^n (z^j g_n(z))\|_R \leq \lambda^n(m+1) C^n B(j+n) \alpha^n (n\beta + R)^{j+n-1}.
\end{align}
To show (\ref{rlimit}), we rewrite the entire expression in (\ref{rlimit}) using $g_n(z)$ and estimate that
\begin{align}\label{Tneqn} \notag
    & \left \| T^n \left( \sum_{j=\floor*{ n^{1.2}} + 1}^\infty \frac{(rz)^j}{j!} g_n(z) \right) \right\|_R \\ \notag
    & \leq   \sum_{j=\floor*{ n^{1.2}} + 1}^\infty \frac{|r|^j}{j!} \left\| T^n (z^j g_n(z))\right\|_R  \\
    & \leq  \sum_{j=\floor*{ n^{1.2}} + 1}^\infty \frac{|r|^j}{j!} \lambda^n (m+1)C^n B (j+n) \alpha^n (n\beta + R)^{j+n-1} \text{ by (\ref{Tnzj}).}
\end{align}
The ratio of successive  terms in the above sum, after simplifying, is
\begin{align}\label{ratio}
    \frac{|r|(j+n+1)(n\beta + R)}{(j+n)(j+1)}.
\end{align}
Since $j\geq \floor*{ n^{1.2}} + 1$, the expression in (\ref{ratio}) is less than $1/2$ for large enough $n$.  The leading term in the sum in (\ref{Tneqn}) is

\begin{align*}
   & \frac{|r|^{\floor*{n^{1.2}} + 1}}{(\floor*{ n^{1.2}} + 1)!} (m+1)C^n B (\floor*{ n^{1.2}} + 1+n) \alpha^n
    (n\beta + R)^{\floor*{ n^{1.2}} + n} \\
    & \leq |r|\frac{|r|^ {n^{1.2}}}{( n^{1.2})!} (m+1)C^n B(3n^{1.2})\alpha^n (2n\beta)^{n^{1.2}+n} \text{ whenever $n\beta >R$}\\
    & \leq |r|\frac{|re|^{n^{1.2}}}{n^{1.2n^{1.2}}} (m+1)C^n B(3n^{1.2})\alpha^n (2n\beta)^{n^{1.2}+n} \\
    & = \frac{1}{n^{0.2n^{1.2}-n-1.2}} |r| (m+1) 3B (2C\alpha \beta)^n |2re\beta|^{n^{1.2}}\\
    & <  \frac{1}{n^{0.2n^{1.2}-n-1.2}} |12r^2(m+1)BC\alpha\beta^2e|^{n^{1.2}}  \\ 
    & < \frac{1}{n^{0.1n^{1.2}}} |12r^2(m+1)BC\alpha\beta^2e|^{n^{1.2}} \text{ for large enough $n$} \\
     & = \left| \frac{12r^2(m+1)BC\alpha\beta^2 e}{n^{0.1}}\right|^{n^{1.2}},
\end{align*}
which converges to 0 as $n \to \infty$. Thus the sum in (\ref{Tneqn}) converges to 0 by the ratio test, which finishes the proof.
\end{proof}

%%%%%%%%%%%%%%%%%%%%%%%%%%
%%%%% e TAYLOR LEMMA %%%%%
%%%%%%%%%%%%%%%%%%%%%%%%%%

In the proof of the previous lemma, Taylor polynomials of $e^{-rz}$ and $e^{rz}$ have played an important role.  Our next lemma also addresses these Taylor polynomials to aid us in our construction of a hypercyclic function.

\begin{lemma}\label{eTaylorLemma}
Let $N$ and $M$ be positive integers with $M<N$.  Let $r$ be a nonzero complex number and let $R>0$ and $0<\sigma<1$ such that $R<M^{1-\sigma}$.  If $M^{0.5\sigma} > 4|r|e$, then
\[ \left\| 1-\sum_{i=0}^M \frac{(-rz)^i}{i!} \sum_{j=0}^N \frac{(rz)^j}{j!}\right\|_R < \frac{e^M}{M^{0.5\sigma M}} .\]
\end{lemma}

\begin{proof}
Let 
\begin{eqnarray}\label{akDef}
\sum_{i=0}^M \frac{(-z)^i}{i!} \sum_{j=0}^N \frac{z^j}{j!} = \sum_{k=0}^{N+M} a_k z^k.
\end{eqnarray}
We first obtain a bound on the coefficients $a_k$. It is clear that $a_0=1$.  By a change of variable $k=j+i$, rewriting the double sum in terms of $z^k$ we have that
\begin{align*}
    \sum_{i=0}^M \left(\frac{(-z)^i}{i!} \sum_{j=0}^N \frac{z^j}{j!}\right) = \sum_{k=0}^{N+M}  \sum_{i=\max(0,k-N)}^{\min(k,M)} \frac{(-1)^i}{i!(k-i)!} z^k.
\end{align*}
Thus for $k\geq 1$,
\begin{align}
    a_k = \sum_{i=\max(0,k-N)}^{\min(k,M)} \frac{(-1)^i}{i!(k-i)!}.\label{akEqn}
\end{align}
We now compute $a_k$ for three different ranges of the values of $k$. First, when $1\leq k\leq M$, we have that 
\begin{align*}
    a_k = \sum_{i=0}^{k} \frac{(-1)^i}{i!(k-i)!} = \frac{1}{k!} \sum_{i=0}^{k} (-1)^i\binom{k}{i} = 0,
\end{align*}
by the Binomial Theorem.

Second, when $M<k\leq N$, we use the formula for a truncated alternating sum of binomial coefficients to see that, in this case, equation (\ref{akEqn}) becomes
\begin{align}
    a_k = \sum_{i=0}^{M} \frac{(-1)^i}{i!(k-i)!} = \frac{1}{k!} \sum_{i=0}^{M} (-1)^i\binom{k}{i} = (-1)^M \frac{1}{k!} \binom{k-1}{M}.\label{akEqn2}
\end{align}

Lastly, when $N+1\leq k \leq N+M$, by using the same formula to derive (\ref{akEqn2}), equation (\ref{akEqn}) becomes 
\begin{align} \notag
    a_k = \sum_{i=k-N}^{M} \frac{(-1)^i}{i!(k-i)!} & = \sum_{i=0}^{M} \frac{(-1)^i}{i!(k-i)!} - \sum_{i=0}^{k-N-1} \frac{(-1)^i}{i!(k-i)!}\\\notag
    & = \frac{1}{k!} \sum_{i=0}^{M} (-1)^i\binom{k}{i} - \frac{1}{k!}\sum_{i=0}^{k-N-1} (-1)^i\binom{k}{i} \\
    & = (-1)^M \frac{1}{k!}\binom{k-1}{M} - (-1)^{k-N}\frac{1}{k!}\binom{k-1}{k-N-1}.\label{akEqn3}
\end{align}
If $M\leq  \floor{\frac{k}{2}}$, then since $k-N\leq M$, we have that $\binom{k-1}{k-N-1}\leq \binom{k-1}{M}$.  In the other case that $M> \floor{\frac{k}{2}}$, then since $N\geq M$, we also have that $\binom{k-1}{k-N-1} = \binom{k-1}{N}\leq \binom{k-1}{M}$.  By combining this observation with equations (\ref{akEqn3}) and (\ref{akEqn2}), we have that for $M<k\leq N+M$,
\begin{align}\label{akBound} \notag
    |a_k| \leq \frac{2}{k!} \binom{k-1}{M} & =\frac{2}{k} \frac{1}{M!(k-M-1)!}\\ \notag
    & =\frac{2}{k} \frac{k(k-1)\cdots (k-M)}{M!k!}\\ \notag
    & < \frac{2k^M}{M!k!}\\
    & <\frac{2k^M e^{M+k}}{M^M k^k},  \text{ because $\displaystyle\frac{n^n}{n!}<e^n$ when $n>0$.}
\end{align}

Now suppose that $r$ is a nonzero complex number, that $0<R<M^{1-\sigma}$ where $0<\sigma<1$, and that $M^{0.5\sigma} > 4|r|e$. We now use (\ref{akDef}) to compute that

\begin{align*}
     \left\| 1-\sum_{i=0}^M \frac{(-rz)^i}{i!} \sum_{j=0}^N \frac{(rz)^j}{j!}\right\|_R & = \left\| \sum_{k=M+1}^{M+N} a_k r^k z^k \right\|_R \\
     & \leq \sum_{k=M+1}^{M+N} |a_k| |r|^k R^k \\
     & < \sum_{k=M+1}^{M+N} \frac{2k^M e^{M+k}}{M^M k^k} |r|^k \frac{M^k}{M^{\sigma k}} \text{ by (\ref{akBound})} \\
     & = e^M \sum_{k=M+1}^{M+N} \frac{2|r|^k e^k }{M^{\sigma k}} \left(\frac{k^M M^k}{k^kM^M}\right) \\
     & \leq \frac{e^M}{M^{0.5\sigma M}} \sum_{k=M+1}^{M+N} \frac{2|r|^k e^k }{M^{0.5\sigma k}}.
\end{align*}
Since $\ds \frac{|r|e}{M^{0.5\sigma}}< \frac{1}{4}$, the geometric sum in the above expression has a leading term less than $1/2$ and a common ratio less than $1/4$ .  Hence the sum is no greater than one, and we have our desired inequality.
\end{proof}

%%%%%%%%%%%%%%%%%%%%%%%%%%%%
%%%%% QUOTIENT POLYS %%%%%%%
%%%%%%%%%%%%%%%%%%%%%%%%%%%%

\section{Quotient Polynomials}\label{QuotientPolys}

It turns out that quotient polynomials play an important role in our proof of our main result Theorem \ref{HypThm}. In this section, we use the lemmas in the previous section to obtain a few preliminary results on quotient polynomials. Our first lemma deals with the remainder polynomial that is obtained from our use of the division algorithm. 

\begin{lemma}\label{RemainderLemma}
Suppose $f$ is a polynomial with simple zeros $\alpha_1,\ldots,\alpha_m$ contained in a disk $|z|\leq R$.  There exists a constant $\omega$, independent of $R$, with the following property: if $g_n$, $q_n$, and $r_n$ are polynomials for which $g_n=fq_n+r_n$ and $r_n=r_{0,n}+r_{1,n} z +\cdots+r_{m-1,n} z^{m-1}$, then
\[ \max\{|r_{i,n}| : 0\leq i\leq m-1\} \leq \omega \|g_n\|_R.\]
Consequently, if $R_n > 0$ with $\lim_{n\to \infty} R_n=\infty$ and
\[ \lim_{n\to\infty} R_n^m \|g_n\|_{R_n} = 0,\]
then $\lim_{n\to\infty}  \|r_n\|_{R_n} = 0$ and $\lim_{n\to\infty}  \|q_n\|_{R_n} = 0$.

\end{lemma}

\begin{proof}
Let $V$ be the Vandermonde matrix
\[ V= 
\begin{bmatrix}
1 & \alpha_1 & \dots & \alpha_1^{m-1}\\
1 & \alpha_2 & \dots & \alpha_2^{m-1}\\
\vdots & \vdots & \ddots & \vdots\\
1 & \alpha_m & \dots & \alpha_m^{m-1}\\
\end{bmatrix}.
\]
Let $\overrightarrow{r_n}=\begin{bmatrix}
r_{0,n} & r_{1,n}  & \dots & r_{m-1,n}
\end{bmatrix}^T$
and $\overrightarrow{g_n} = \begin{bmatrix}
g_n(\alpha_1) & g_n(\alpha_2) & \dots & g_n(\alpha_{m})
\end{bmatrix}^T$. Since $g_n(\alpha_i)=r_n(\alpha_i)$ for $1\leq i \leq m$, we have that $V\overrightarrow{r_n}=\overrightarrow{g_n}$.  Since the $\alpha_i$ are distinct, $V$ is invertible and $\overrightarrow{r_n}=V^{-1}\overrightarrow{g_n}$.  Then by applying the sup-norm and its corresponding induced matrix norm $\| \cdot \|$, we have that 
\[ \|\overrightarrow{r}\|_\infty \leq \|V^{-1}\|\cdot \|\overrightarrow{g_n}\|_\infty.  \]
Since each $\alpha_i$ is in the disk $|z|\leq R$, we have that $|g_n(\alpha_i)|\leq \|g_n\|_R$, which finishes the first part of the proof by setting $\omega=\|V^{-1}\|$.

To finish the proof of our lemma, suppose $R_n > 0$ with $\lim_{n\to \infty} R_n=\infty$ and
$\lim_{n\to\infty} R_n^m \|g_n\|_{R_n} = 0$. Since $R_n\to \infty$, the zeros of $f$ are contained in all but finitely many of the disks $|z|\leq R_n$. Hence $\|\overrightarrow{g_n}\|_\infty \leq \|g_n\|_{R_n}$ for large enough $n$, so eventually
\begin{align}\notag
    \|r_n\|_{R_n} \leq \sum_{i=0}^{m-1} |r_{i,n}| R_n^i \leq m \|V^{-1}\| R_n^m \|g_n\|_{R_n} .
\end{align}
Hence  $\lim_{n\to\infty}\|r_n\|_{R_n}=0$.  

Since $f$ is a polynomial, $\ds\lim_{n\to\infty}\left( \sup_{|z|=R_n} \frac{1}{|f(z)|}\right) = 0$. Then since $q_n=(g_n-r_n)/f$, we have that
\begin{align}\notag
    \lim_{n\to\infty} \|q_n\|_{R_n} \leq \lim_{n\to\infty}\left( \sup_{|z|=R_n} \frac{1}{|f(z)|}\right) ( \|g_n\|_{R_n} + \|r_n\|_{R_n}) = 0.
\end{align}
\end{proof}

Our next proposition is used to determine the polynomials $q_n$ that we use to construct our infinite product in the case that our differential operator $T=\varphi(D)$ is such that $\varphi(0)\not = 0$.

\begin{proposition}\label{Part1}
Suppose $T = \varphi (D)$ is a differential operator for some non-constant entire function $\varphi$ of exponential type with $\varphi (0) \not= 0$. Let $f$ be a polynomial of degree $m$ with simple zeros and $f(0)\not =0$. Let $p$ be a polynomial of degree less than $m$. Then there exist polynomials $q_n$ with $n \leq \deg q_n < n^{1.3}$ and the following properties hold:
\begin{enumerate}[label={\upshape(\alph*)}]
    \item as $n \to \infty, \   q_n\to 0$ uniformly on compact subsets of $\C$, and for large enough $n$, any zero of the polynomial $q_n+1$ is greater in modulus than $n^{0.7}$,
    \item $T^n ((q_n+1)f)\to  p$ as $n \to \infty$,
    \item for large enough $n$, the coefficient $q_{n,1}$ of $z$ in the polynomial $q_n$ satisfies $|q_{n,1}|<1/n$, and
    \item for large enough $n$, $q_n(0)=0$ and the zeros of the polynomial $q_n+1$ are simple.
\end{enumerate}
\end{proposition}

\begin{proof}
Let $S:\mathcal{P}\to\mathcal{P}$ be the right inverse for $T$ given by Lemma \ref{InverseLemma}, and let $r$ be a non-zero complex number given by Lemma \ref{rLemma}. Define $h_n$ to be
\[h_n(z) = S^n p(z) + (f(z) - S^n p(z))  \displaystyle\sum\limits_{i=0}^{n-m-1} \frac{(-rz)^i}{i!} \displaystyle\sum\limits_{j=0}^{\floor*{ n^{1.2}}}  \frac{(rz)^j}{j!}. \]
Let $q_n$ be the quotient polynomial of the division of $h_n - f$ by $f$, which means that $h_n -f = f q_n +r_n$ and $\deg r_n < m$. We now show that $q_n$ satisfies the conditions of the proposition.  Since the degree of $h_n$ is $\floor*{ n^{1.2}}  + n-1$, the degree of $q_n$ is $\floor*{ n^{1.2}}  + n-m-1$, which is greater than $n$, and is less than $n^{1.3}$ for large enough $n$.

We now estimate the norm $\|h_n-f\|_{n^{0.7}}$. First we observe
\begin{align}\label{hn-fpreview} \notag
    \| h_n -f \|_{n^{0.7}} & = \left \| (S^n p(z) - f(z)) \left( 1 -  \displaystyle\sum\limits_{i=0}^{n-m-1} \frac{(-rz)^i}{i!}  \displaystyle\sum\limits_{j=0}^{\floor*{n^{1.2}}} \frac{(rz)^j}{j!} \right)\right\|_{n^{0.7}}\\ 
    & \leq \| S^n p(z) - f(z) \|_{n^{0.7}} \left \| 1 -  \displaystyle\sum\limits_{i=0}^{n-m-1} \frac{(-rz)^i}{i!}  \displaystyle\sum\limits_{j=0}^{\floor*{n^{1.2}}} \frac{(rz)^j}{j!} \right\|_{n^{0.7}}.
\end{align}
By Lemma \ref{InverseLemma}, there is a constant $C>1$ for which the coefficients of the polynomial $S^n p(z) - f(z)$ are bounded above in modulus by $C^n$. Hence
\begin{align}\label{Snp-f bound}
    \| S^n p(z) - f(z)\|_{n^{0.7}} \leq (m+1) C^n n^{0.7m}.
\end{align}
For large enough $n$, we have $n^{0.7}<(n-m-1)^{1-0.2}$, so we can apply (\ref{Snp-f bound}) and Lemma \ref{eTaylorLemma} on (\ref{hn-fpreview}) with $\sigma = 0.2, M=n-m-1, N=\floor*{n^{1.2}}$, and $R=n^{0.7}$ to obtain that
\begin{align}\label{hn-f}
     \| h_n -f \|_{n^{0.7}} \leq (m+1) C^n n^{0.7m} \cdot \frac{e^{n-m-1}}{(n-m-1)^{0.1(n-m-1)}}.
\end{align}
Consequently,
\begin{align*}
    \lim_{n\to\infty} n^{0.7m}  \| h_n - f\|_{n^{0.7}} \leq \lim_{n\to\infty} \frac{(m+1) C^n n^{1.4m}e^{n-m-1}}{(n-m-1)^{0.1(n-m-1)}} =0.
\end{align*}
Thus by applying Lemma \ref{RemainderLemma} with $R_n=n^{0.7}$ and $g_n=h_n-f$, both $\|r_n\|_{n^{0.7}}\to 0$ and $\|q_n\|_{n^{0.7}}\to 0$ as $n\to \infty$. Whenever $\|q_n\|_{n^{0.7}} < 1$ we have
\begin{align*}
    \inf_{|z|\leq n^{0.7}} |1+q_n(z)|  \geq \inf_{|z| \leq n^{0.7}} \big| 1 - |q_n(z)| \big| \geq 1 - \|q_n\|_{n^{0.7}}.
\end{align*}
Thus for large enough $n$, the polynomial $q_n + 1$ has no zeros in the disk $|z|\leq n^{0.7}$. This finishes the proof of property (a).

To prove property (b), let $R>0$ be given.  Since $T$ is linear and  $T^n S^n p(z) = p(z)$, Lemma \ref{rLemma} implies that $T^n h_n \to p$.  By Lemma \ref{RemainderLemma}, there exists a constant $\omega$ for which the coefficients of $r_n$ are bounded above in modulus by $\omega \|h_n - f\|_{n^{0.7}}$.  Then by Lemma \ref{PolyGrowthLemma}, there are positive constants $\alpha$ and $\beta$ such that 
\begin{align*}
    \|T^n(r_n)\|_R & \leq m \omega \|h_n - f\|_{n^{0.7}} \alpha^n (n\beta+R)^{m-1} \\
    & < m \omega \frac{(m+1) C^n n^{0.7m}e^{n-m-1}}{(n-m-1)^{0.1(n-m-1)}} \alpha^n (n\beta+R)^{m-1}  \text{ by (\ref{hn-f}) }\\
    & \to 0 \text{ as } n \to \infty. &
\end{align*}
Thus $T^n(r_n)\to 0$ as $n \to \infty$.  Since $(q_n+1)f=h_n-r_n$, we have that
\begin{align*}
    T^n ((q_n+1)f)  = T^n(h_n)-T^n(r_n) \to p.
\end{align*}
This proves property (b).

For property (c), we use the fact that the quotient polynomial $q_n$ is the principal part of $\ds\frac{h_n-f}{f}$ at $\infty$; see \cite[page 180]{Gamelin}. Thus by the formula for Laurent series coefficients \cite[Equation VI.1.4]{Gamelin}, if the zeros of $f$ are contained within the open disk $|z|<R$, then the coefficient $q_{n,1}$ of $z$ in $q_n$ satisfies
\begin{align*}
    |q_{n,1}| = \left| \frac{1}{2\pi i}\oint_{|z|=R} \frac{h_n(z) -f(z)}{f(z)} \frac{1}{z^2}\, dz \right | \leq \frac{\|h_n-f\|_R}{\inf\limits_{|z|=R} |f(z)|} \frac{1}{R}.
\end{align*}
By setting $R=n^{0.7}$ and applying inequality (\ref{hn-f}), we have
\begin{align*}
    \lim_{n\to\infty} n|q_{n,1}| \leq \lim_{n\to\infty} \frac{n}{\inf\limits_{|z|=n^{0.7}} |f(z)|} \frac{(m+1) C^n n^{0.7m}e^{n-m-1}}{(n-m-1)^{0.1(n-m-1)}} \frac{1}{n^{0.7}} = 0.
\end{align*}
This finishes the proof of property (c).

For property (d), we estimate $|q_n(0)|$ using a similar argument above as
\begin{align}\label{qn0} 
    |q_n(0)| = \left| \frac{1}{2\pi i}\oint_{|z|=R} \frac{h_n(z) -f(z)}{f(z)} \frac{1}{z}\, dz \right | \leq \frac{\|h_n-f\|_R}{\inf\limits_{|z|=R} |f(z)|}.
\end{align}
Now if $\mu$ is the maximum modulus of the coefficients of $f$, then Lemma $\ref{PolyGrowthLemma}$ implies that there exist positive constants $\alpha$ and $\beta$ such that $\|T^n f\|_R \leq \mu(m+1)\alpha^n (n\beta+R)^m$.  Hence
\begin{align}\label{(d)eqn} \notag
    \|T^n(q_n(0)f)\|_R
    & \leq |q_n(0)| \mu(m+1)\alpha^n (n\beta+R)^m. \\ 
    & \leq  \frac{\|h_n-f\|_R}{\inf\limits_{|z|=R} |f(z)|} \mu(m+1)\alpha^n (n\beta+R)^m \text{ by (\ref{qn0}).}
\end{align}
The maximum modulus principle implies that eventually $\|h_n -f \|_R$ $\leq$ $\|h_n-f\|_{n^{0.7}}$. Thus by applying inequality (\ref{hn-f}) to inequality (\ref{(d)eqn}), one can show that $\|T^n(q_n(0)f)\|_R$ converges to 0. Combining this result with property (b), we have that
\begin{align*}
    T^n((q_n-q_n(0)+1)f) \to p.
\end{align*}
Thus property (b) holds if we replace $q_n$ with $q_n-q_n(0)$.  Since $\|q_n-q_n(0)\|_{n^{0.7}} \leq 2\|q_n\|_{n^{0.7}}$, property (a) also holds by making this same replacement.  Since property (c) does not depend on the constant term of $q_n$, we may replace $q_n$ by $q_n-q_n(0)$ in properties (a), (b) and (c) and have them remain true.  If the zeros of $q_n-q_n(0)+1$ are simple, we are finished, so suppose they are not.  By using small enough perturbations of any repeated zeros, properties (a) and (c) still hold.  Since $T$ is continuous, by using small enough perturbations of any repeated zeros, property (b) still holds as well.  Thus by replacing $q_n$ with $q_n-q_n(0)$ and perturbing the zeros if necessary, we may assume that $q_n(0)=0$ and the zeros of $q_n+1$ are simple.
\end{proof}

We now move on to the other case when $T$ is a differential operator $T=\varphi(D)$ for some non-constant entire function $\varphi$ of exponential type with $\varphi(0)=0$.  In this case, $\varphi(z)=z^k\psi(z)$ for some positive integer $k$ and entire function $\psi(z)$ with $\psi(0)\not = 0$.  This makes $T=\psi(D) D^k$, and the presence of a factor of $D$ makes our proof for the following proposition considerably easier than its analogue Proposition \ref{Part1}.

\begin{proposition}\label{Part2}
Suppose $T = \varphi (D)$ is a differential operator for some non-constant entire function $\varphi$ with $\varphi (0) = 0$.   Let $f$ be a polynomial of degree $m$ with simple zeros and $f(0)\not =0$. Let $p$ be a polynomial of degree less than $m$. Then there exist polynomials $q_n$ with $n \leq \deg q_n < n^{1.3}$ and the following properties hold:
\begin{enumerate}[label={\upshape(\alph*)}]
    \item as $n \to \infty, \   q_n\to 0$ uniformly on compact subsets of $\C$, and for large enough $n$, any zero of the polynomial $q_n+1$ is greater in modulus than $n^{0.7}$,
    \item $T^n ((q_n+1)f)\to  p$ as $n \to \infty$,
    \item for large enough $n$, the coefficient $q_{n,1}$ of $z$ in the polynomial $q_n$ satisfies $|q_{n,1}|<1/n$, and
    \item for large enough $n$, $q_n(0)=0$ and the zeros of the polynomial $q_n+1$ are simple.
\end{enumerate}
\end{proposition}

\begin{proof}
Write $\varphi(z) = z^k \psi(z)$, where $k\geq 1$ and $\psi(0)\not = 0$. Then the operator $T$ can be expressed as $T=\Tilde{T}D^k$, where $\Tilde{T}=\psi(D)$.  Note that $Tg=0$ for any polynomial $g$ of degree at most $k-1$.  Let $\mathcal{P}$ be the collection of polynomials in $H(\C)$ and let $g\in \mathcal{P}$ be given.  By Lemma \ref{InverseLemma}, there is a mapping $\Tilde{S}:\mathcal{P}\to \mathcal{P}$ such that $\Tilde{T}\Tilde{S}g=g$.  Define $A:\mathcal{P}\to \mathcal{P}$ by $Ag(z) = \int_0^z g(w)\, dw$ and $S_n: \mathcal{P}\to \mathcal{P}$ by
\begin{align}\label{Sinv}
    S_n g = A^{kn}\Tilde{S}^n g.
\end{align}
Then 
\[ T^n S_n g = \Tilde{T}^n D^{kn} A^{kn}\Tilde{S}^n g = \Tilde{T}^n \Tilde{S}^n g = g,\]
so $S_n$ is a right inverse for $T^n$ on the set $\mathcal{P}$.

Let $\Tilde{S}^n p = c_{n,0}+\cdots + c_{n,m} z^m$ as in Lemma \ref{InverseLemma}.  Then there is some constant $C>1$ for which
\begin{align}\notag
    \| S_n p\|_{n^{0.7}} & = \left \| \sum_{j=0}^m A^{kn} (c_{n,j} z^j) \right\|_{n^{0.7}} & \\ \notag
    & \leq \sum_{j=0}^m \left\|C^n \frac{z^{j+kn}j!}{(j+kn)!}\right\|_{n^{0.7}}  \text{ by Lemma \ref{InverseLemma} }\\ \notag
    & < \sum_{j=0}^m \frac{C^n n^{0.7(j+kn)}}{(kn)!} &\\ \notag
    & < \frac{(m+1)C^n n^{0.7(m+kn)}}{1} \frac{e^{kn}}{(kn)^{kn}} &\\
    & = \frac{(m+1) C^n n^{0.7m} }{n^{0.3kn}} \frac{e^{kn}}{k^{kn}}. & \label{SnpEst}
\end{align}
Suppose $n$ is large enough so that $kn > m-\deg p$.  Let $q_n$ and $r_n$ be the quotient and remainder of the division of $S_n p$ by $f$. That is, $S_n p = fq_n + r_n$ and $\deg r_n < m$.  Since the degree of $S_n p$ is $kn+\deg p$ by (\ref{Sinv}), it must be that $\deg q_n =kn+\deg p-m$, which is less then $n^{1.3}$ if $n$ is large enough.  Furthermore, since by (\ref{SnpEst})
\begin{align*}
    \lim_{n\to\infty} n^{0.7m} \|S_n p\|_{n^{0.7}} \leq \lim_{n\to\infty} \frac{(m+1) C^n n^{1.4m} }{n^{0.3kn}} \frac{e^{kn}}{k^{kn}} = 0,
\end{align*}
by setting $R_n=n^{0.7}$ and $g_n=S_n p$, Lemma \ref{RemainderLemma} implies that $\lim_{n\to\infty} \|q_n\|_{n^{0.7}} = 0$.  Whenever $\|q_n\|_{n^{0.7}} < 1$ we have
\begin{align*}
    \inf_{|z|\leq n^{0.7}} |1+q_n(z)|  \geq \inf_{|z| \leq n^{0.7}} \big| 1 - |q_n(z)| \big| \geq 1 - \|q_n\|_{n^{0.7}}.
\end{align*}
Thus for large enough $n$, the polynomial $q_n + 1$ has no zeros in the disk $|z|\leq n^{0.7}$. This finishes the proof of property (a).

To prove property (b), we use the fact that whenever $n>m$, $T^n f =0=T^n(r_n)$.  Hence for $n>m$,
\[T^n ((q_n+1)f)=T^n(fq_n)=T^n(fq_n + r_n) =T^n (S_n p)=p,\]
which proves property (b).

For property (c), we replicate the argument in Proposition \ref{Part1}. If the zeros of $f$ are contained in the disk $|z|<R$, then the coefficient $q_{n,1}$ of $z$ in $q_n$ satisfies
\begin{align*}
    |q_{n,1}| = \left| \frac{1}{2\pi i}\oint_{|z|=R} \frac{S_n p(z)}{f(z)} \frac{1}{z^2}\, dz \right | \leq \frac{\|S_n p(z)\|_R}{\inf\limits_{|z|=R} |f(z)|} \frac{1}{R}.
\end{align*}
By setting $R=n^{0.7}$ and applying inequality (\ref{SnpEst}), we have
\begin{align*}
    \lim_{n\to\infty} n|q_{n,1}| \leq \lim_{n\to\infty} \frac{n}{\inf\limits_{|z|=n^{0.7}} |f(z)|} \frac{(m+1) C^n n^{0.7m} }{n^{0.3kn}} \frac{e^{kn}}{k^{kn}} \frac{1}{n^{0.7}} = 0.
\end{align*}
This finishes the proof of property (c).

We now prove property (d).  Since $T^n(fq_n(0))=0$ for $n>m$, property (b) still holds if we replace $q_n$ with $q_n-q_n(0)$.  Since $\|q_n-q_n(0)\|_{n^{0.7}} \leq 2\|q_n\|_{n^{0.7}}$, property (a) also holds by making this same replacement.  Since property (c) does not depend on the constant term of $q_n$, we may replace $q_n$ by $q_n-q_n(0)$ in properties (a), (b) and (c) and have them remain true.  If the zeros of $q_n-q_n(0)+1$ are not simple, then by using small enough perturbations of the zeros, we can assume that properties (a), (b) and (c) still hold for polynomials $q_n$ with $q_n(0)=0$ and such that $q_n+1$ has simple zeros.
\end{proof}

We combine the previous propositions to prove the following proposition which helps us factor our infinite product into linear factors.

\begin{proposition}\label{EndGame}
Suppose $T = \varphi (D)$ is a differential operator for some non-constant entire function $\varphi$.   Let $f$ be a polynomial of degree $m$ with simple zeros and $f(0)\not = 0$. Let $p$ be a polynomial of degree less than $m$. Then there exist polynomials $q_n$ with $n \leq \deg q_n < n^{1.3}$ such that properties (a)-(d) in Proposition \ref{Part2} hold, and the following additional property holds as well:
\begin{enumerate}[label={\upshape(\alph*)}]
    \setcounter{enumi}{4}
    \item Let $a_{n,1},\ldots, a_{n,\deg q_n}$ be the zeros of $q_n+1$. For any positive constants $\epsilon$ and $R$, for sufficiently large $n$ there is an ordering of the zeros such that for any $J$ satisfying $1\leq J \leq \deg q_n$,
\begin{align*}
    \left\|  1 - \prod\limits_{j=1}^{J} \left( 1 - \frac{z}{a_{n,j}} \right)  \right\|_R < \epsilon.
\end{align*}

\end{enumerate}
\end{proposition}

\begin{proof}
Let $R>0$ and $\epsilon > 0$ be given.  Let $\Log z$ be the principal branch of the logarithm defined on $\C\setminus (-\infty,0]$. We may suppose by property (a) that $n$ is large enough so that the zeros $a_{n,1},\ldots, a_{n,\deg q_n}$ of the polynomial $q_n + 1$ satisfy $|a_{n,j}| > R$.  Hence the real part of $\ds 1 -\frac{z}{a_{n,j}}$ is positive on the disk $|z|\leq R$, so the function $\ds\Log \left(1 -\frac{z}{a_{n,j}}\right)$ is analytic on the disk $|z|\leq R$. To prove property (e), it is sufficient to show that for sufficiently large $n$,
\begin{align}
\max_{|z|\leq R} \left| \Log  \prod\limits_{j=1}^{J} \left( 1 - \frac{z}{a_{n,j}}  \right) \right| < \epsilon \label{log switch}
\end{align}
for some ordering of the zeros $a_{n,1},a_{n,2},\dots,a_{n,\deg q_n}$ of $q_n+1$.

We use Taylor series to express
\begin{align}\notag
\Log \prod\limits_{j=1}^{J} \left( 1 - \frac{z}{a_{n,j}} \right) &=
 \sum\limits_{j=1}^{J} \Log \left( 1 - \frac{z}{a_{n,j}} \right) \\ \notag
 &= \sum\limits_{j=1}^{J} \sum\limits_{k=1}^{\infty} \frac{-(z)^k}{k(a_{n,j})^k} \\ \notag
 & = -\sum\limits_{k=1}^{\infty} \frac{z^k}{k} \sum\limits_{j=1}^{J} \left( \frac{1}{(a_{n,j})^k} \right) \\
 &=\left( -z \sum\limits_{j=1}^{J} \frac{1}{a_{n,j}} \right) + \left( -\sum\limits_{k=2}^{\infty} \frac{z^k}{k} \sum\limits_{j=1}^{J} \left( \frac{1}{(a_{n,j})^k} \right) \right). \label{two terms}
\end{align}

Let $q_{n,1}$ be the coefficient of $z$ in $q_n$. Suppose $n$ is so large that we have both
\begin{align}
    |1/a_{n,j}|<1/n^{0.7} \text{ \hskip .2in and \hskip .2in } |q_{n,1}| <1/n^{0.7}, \label{anjEst}
\end{align}
by properties (a) and (c).  Since $q_n(0)=0$ by property (d), $q_n(z) + 1=\ds \prod_{j=1}^{\deg q_n}\left(1-\frac{z}{a_{n,j}}\right)$ and $q_{n,1} = -\ds \sum_{j=1}^{\deg q_n} \frac{1}{a_{n,j}}$.  Now as a consequence of the Polygonal Confinement Theorem \cite[Lemma 3.1]{Rosenthal}, the inequalities in (\ref{anjEst}) imply that we can arrange the $a_{n,j}$ so that if $1\leq J \leq \deg q_n$, then
\begin{align}
\left| \sum\limits_{j=1}^{J} \frac{1}{a_{n,j}} \right|<\frac{\sqrt{5}}{n^{0.7}}.\label{PolyConfine}
\end{align}

Therefore, 

\begin{align*}
    & \max_{|z|\leq R} \left| \Log \left( \prod\limits_{j=1}^{J} \left( 1 - \frac{z}{a_{n,j}} \right) \right) \right|  \\
    & \leq\max_{|z|\leq R} \left| -z \sum\limits_{j=1}^{J} \frac{1}{a_{n,j}} \right| + \max_{|z|\leq R} \left|- \sum\limits_{k=2}^{\infty} \frac{ z^k}{k} \sum\limits_{j=1}^{J} \left( \frac{1}{(a_{n,j})^k} \right) \right|  \text{ by (\ref{two terms}) } \\
    & \leq
 \frac{\sqrt{5}R}{n^{0.7}} + \sum\limits_{k=2}^{\infty} \frac{R^k}{k} \sum\limits_{j=1}^{J} \frac{1}{|a_{n,j}|^k}  \text{ by (\ref{PolyConfine}) }\\
& < \frac{\sqrt{5}R}{n^{0.7}} + \sum\limits_{k=2}^{\infty} R^k n^{1.3}  \frac{1}{n^{0.7k}} \text{ since $J<n^{1.3}$ and by (\ref{anjEst})} \\
& \leq \frac{\sqrt{5}R}{n^{0.7}} + \frac{R^2}{n^{0.1}} \sum\limits_{k=0}^{\infty} \left( \frac{R }{n^{0.7}} \right)^k.
\end{align*}
The sum in the previous line is a convergent geometric series, and the other terms converge to 0 as $n \to \infty$. Therefore (\ref{log switch}) is fulfilled for sufficiently large $n$, which completes the proof.
\end{proof}

%%%%%%%%%%%%%%%%%%%%%%%
%%%%% MAIN RESULT %%%%%
%%%%%%%%%%%%%%%%%%%%%%%
\section{Main Result}\label{MainResult}

We can now use Proposition \ref{EndGame} to prove our main result.  We borrow the multiplicative techniques from \cite{ChanWalmsley} to do our construction of a particular hypercyclic function as an infinite product.

\begin{theorem}\label{HypThm}
Suppose $T=\varphi(D): H(\C)\to H(\C)$ is a differential operator for some non-constant entire function $\varphi$ of exponential type.  Then there exists a function $f$ in $H(\C)$ which is hypercyclic for $T$ such that
\[ f(z)=\prod_{j=1}^\infty \left(1-\frac{z}{a_j}\right), \]
where $a_j$ are nonzero complex numbers.
\end{theorem}

\begin{proof}
We use an inductive process to choose polynomials $\{q_j\}$ for which the function $f=\prod_{j=1}^\infty \left( q_j +1\right)$ is hypercyclic for $T$. 

Let $T=\sum_{j=0}^\infty a_j D^j$ and $J=\min\{j\in \mathbb{N}: a_j\not =0\}$. Let $p_1,p_2,p_3,\ldots$ be a dense sequence of non-zero polynomials such that $p_1=a_J$, $p_2=1$, and $\deg p_j \leq j-1$ for all $j>2$. We first discuss the case when $J=0$. Let $b=1/(|a_0|+|a_1|)$ and define $q_1(z)=bz$.  Then $q_1(0)=0$, $\deg q_1 > \deg p_2$, and $T(q_1(z) + 1) = a_0(bz+1) + a_1b$.  Therefore
\begin{align*}
    \| T(q_1(z)+1) - p_1(z)\|_1 = \| a_0 bz + a_1 b\|_1 = b \| a_0 z+a_1\|_1 \leq 1.
\end{align*}
Secondly for the case when $J\not = 0$, define $q_1(z)=z^J/J!$. Then $q_1(0)=0$, $\deg q_1 > \deg p_2$, and 
\begin{align*}
    T(q_1(z)+1)=a_J D^J (z^J/J! + 1) = a_J=p_1(z).
\end{align*}
In either case, $q_1$ is a polynomial of degree greater than $\deg p_2$, $q_1(0)=0$, $q_1(z)+1$ has simple zeros, and 
\[ \|T(q_1 + 1) -p_1\|_1 \leq 1.\]

Let $n_1=1$, and inductively suppose for $k\geq 2$ we have found increasing integers $n_1<n_2<\dots < n_{k-1}$ and polynomials $q_1,q_2,\ldots, q_{k-1}$.  We now define the integer $n_k$ and polynomial $q_k$ as follows.

Let $f_{k-1}=\prod_{j=1}^{k-1} \left( q_j + 1\right)$.  By the continuity of $T^{n_1}, T^{n_2},\ldots, T^{n_{k-1}}$ at $f_{k-1}$, there exists a $\delta > 0$ such that for all integers $j$ with $1\leq j \leq k-1$, and all entire functions $g$,
\begin{align}
    \| g-1\|_k < \delta \implies \| T^{n_j} (f_{k-1}g) - T^{n_j} (f_{k-1})\|_k < 2^{-k}.\label{delta}
\end{align}
By applying Proposition \ref{EndGame} with $f=f_{k-1}$ and $p=p_k$, there exist an integer $n_k > n_{k-1}$ and a polynomial $q_k$ of degree greater than $k$ such that
\begin{enumerate}[label={\upshape(\alph*)}]
    \item $\|q_k\|_k < \min\{\delta,k^{-2}\}$, and each zero of $q_k +1$ is greater in modulus than the maximum modulus of any zero of $q_{k-1}+1$,
    \item $\|T^{n_k}((q_k+1)f_{k-1}) - p_k \|_k < \ds 2^{-k}$,
    \item the modulus of the coefficient of $z$ in $q_k$ is less than $1/k$,
    \item $q_k(0)=0$ and $q_k+1$ has simple zeros, and
    \item if $a_{k,1},\ldots, a_{k, \deg q_k}$ are the zeros of $q_k+1$, there is an ordering of the zeros of $q_k + 1$ such that for any $j$ satisfying $1\leq j \leq \deg q_k$,
\begin{align*}
    \left\| 1 - \prod\limits_{i=1}^{j} \left( 1 - \frac{z}{a_{k,i}} \right) \right\|_k < \frac{1}{\|f_{k-1}\|_{k-1}}\frac{1}{2^{k-1}}.
\end{align*}
\end{enumerate}

Now let
\begin{align}\label{feqn}
    f(z)=\displaystyle\prod_{j=1}^\infty \left (q_j(z)+1\right).
\end{align}
By \cite[Theorem 5.9]{Conway}, property (a) implies that the infinite product in (\ref{feqn}) converges uniformly on any compact disk $|z|\leq k$ since 
\begin{align*}
    \sum_{j=1}^\infty \|q_j\|_k & = \sum_{j=1}^{k-1} \|q_j\|_k + \sum_{j=k}^\infty \|q_j\|_k \\
    & \leq \sum_{j=1}^{k-1} \|q_j\|_k + \sum_{j=k}^\infty \|q_j\|_j\\
    & \leq  \sum_{j=1}^{k-1} \|q_j\|_k + \sum_{j=k}^\infty \frac{1}{j^2}\\
    & < \infty.
\end{align*}
Thus $f$ is an entire function.  We show that $f$ is hypercyclic for $T$.  

Let 
\begin{align}\label{fkeqn}
    f_{k}=\prod_{j=1}^{k} \left( q_j + 1\right) = (q_k + 1) f_{k-1}.
\end{align}

Since $\|(q_k+1)-1\|_k < \delta$ by property (a), definition (\ref{fkeqn}) and statement (\ref{delta}) imply that whenever $1\leq j\leq k-1,$
\begin{align}
    \|T^{n_j}f_k - T^{n_j} f_{k-1} \|_k < 2^{-k}.\label{Tnj}
\end{align}
Therefore
\begin{align}\notag
    \|T^{n_k}f-T^{n_k}f_k\|_k &= \|T^{n_k}\left(\lim_{m\to \infty}f_m\right)-T^{n_k}f_k\|_k\\ \notag
    & =  \lim_{m\to\infty} \|T^{n_k}f_m-T^{n_k}f_k\|_k\\ \notag
    & \leq \lim_{m\to\infty} \sum_{j=k}^{m-1} \| T^{n_k}f_{j+1}-T^{n_k}f_j\|_k\\ \notag
    & < \lim_{m\to\infty} \sum_{j=k}^{m-1} \frac{1}{2^{j+1}}  \text{ by (\ref{Tnj}) }\\
    & = \frac{1}{2^k}.\label{Tnkf}
\end{align}
Thus
\begin{align*}
    \|T^{n_k}f-p_k\|_k & \leq \|T^{n_k}f-T^{n_k}f_k\|_k + \|T^{n_k}f_k-p_k\|_k \\
    & \leq \frac{1}{2^k} + \frac{1}{2^k}  \text{ by (\ref{Tnkf}) and property (b).}
\end{align*}
Since the sequence of polynomials $\{p_k\}_{k=1}^\infty$ is dense in $H(\C)$, $f$ is hypercyclic for $T$.  What remains to show is that $f$ can be factored into linear factors. 

Let $d_0=0$ and $d_j=\deg f_j$ for $j\geq 1$, so that by (\ref{fkeqn}) we have $d_j-d_{j-1}= \deg(q_j+1)$ for $j\geq 1$.  Let $a_1, a_2, \ldots$ be the zeros of the polynomials $\{q_j+1\}$ ordered so that if $1\leq j$ and $1+d_{j-1}\leq i \leq d_j$, then we have that
\begin{enumerate}[label={\upshape(\roman*)}]
    \item $a_i$ is a zero of $q_j+1$, and
    \item the zeros of $q_j+1$ are arranged so that property (e) holds.
\end{enumerate}
We now show that for any positive integer $k$,
\begin{align}
    \lim_{j\to \infty} \left\| f-\prod_{i=1}^j \left(1-\frac{z}{a_i}\right) \right\|_k = 0.\label{factor}
\end{align}

Let $\epsilon>0$ be given.  There exists a positive integer $N_1 > k$ such that $2^{-n}<\epsilon/2$ whenever $n\geq N_1$.  Since $f=\lim_{n\to\infty}f_n$, there exists a positive integer $N_2$ such that $\|f-f_n\|_k < \epsilon/2$ whenever $n\geq N_2$.

Let $N=\max\{N_1,N_2\}$ and suppose that $j > \deg f_N$.  Then there is some $n\geq N$ for which $\deg f_n < j \leq \deg f_{n+1}$.  We now factor the product in (\ref{factor}) as
\begin{align*}
    \prod_{i=1}^j \left(1-\frac{z}{a_i}\right) = f_n(z) \prod_{i=1+\deg f_n}^j \left(1-\frac{z}{a_i}\right).
\end{align*}
By (i), the zeros $a_i$ for which $1+\deg f_n \leq i \leq j$ are zeros of $q_{n+1}+1$.  Then by (ii) and property (e), since $n>k$ we have that
\begin{align*}
    \left \| 1-\prod_{i=1+\deg f_n}^j \left(1-\frac{z}{a_i}\right) \right\|_k \leq \left \| 1-\prod_{i=1+\deg f_n}^j \left(1-\frac{z}{a_i}\right) \right\|_n < \frac{1}{\|f_n\|_n} \frac{1}{2^n}.
\end{align*}
Therefore, for $j > \deg f_N$,
\begin{align*}
    \left\| f(z)-\prod_{i=1}^j \left(1-\frac{z}{a_i}\right) \right\|_k & \leq \|f-f_n\|_k + \left\| f_n(z) - f_n(z) \prod_{i=1+\deg f_n}^j \left(1-\frac{z}{a_i}\right)\right\|_k\\
    & \leq \frac{\epsilon}{2} + \|f_n\|_n \left \| 1-\prod_{i=1+\deg f_n}^j \left(1-\frac{z}{a_i}\right) \right\|_n \\
    & \leq \frac{\epsilon}{2} + \|f_n\|_n \frac{1}{\|f_n\|_n}\frac{1}{2^n}\\
    & <\epsilon.
\end{align*}
Thus (\ref{factor}) is verified, and since the positive integer $k$ is arbitrary, we have that $f = \ds \prod_{i=1}^\infty \left(1-\frac{z}{a_i}\right),$ and conclude the proof of the theorem.
\end{proof}

\section{Conclusion}
To conclude our discussion, we point out that we can modify the above proof of Theorem \ref{HypThm} to obtain the following generalization:

{\it For each integer $n\geq 1$, let $T_n=\varphi_n(D):H(\C)\to H(\C)$ be a differential operator for some non-constant entire function $\varphi_n$ of exponential type. Then there exists a function $f$ in $H(\C)$ which is hypercyclic for each $T_n$ such that
\[ f(z)=\prod_{j=1}^\infty \left(1-\frac{z}{a_j}\right), \]
where $a_j$ are nonzero complex numbers.}

To describe the modifications of the proof in Theorem \ref{HypThm}, we replace $T$ with $T_1$ and use the same dense sequence of polynomials $p_1, p_2, p_3, \ldots$ as in the beginning of the proof.  The first step of the induction process remains unchanged.  Let $n_1=1$, and inductively suppose for $k\geq 2$ we have found increasing integers $n_1<n_2<\dots < n_{k-1}$ and polynomials $q_1,q_2,\ldots, q_{k-1}$.  We now define the integer $n_k$ and polynomial $q_k$ as follows.

Let $f_{k-1}=\prod_{j=1}^{k-1} \left( q_j + 1\right)$, let $i$ be an integer satisfying $1\leq i\leq k-1$, and let $j$ be an integer such that $i\leq j \leq k-1$.  By the continuity of each operator $T_{i}^{n_j}$ at $f_{k-1}$, there exists a $\delta > 0$ such that  for all integers $i$ with $1\leq i\leq k-1$, for all integers $j$ with $i\leq j\leq k-1$, and for every entire function $g$,
\begin{align*}
    \| g-1\|_k < \delta \implies \left\| T_{i}^{n_j} (f_{k-1}g) - T_{i}^{n_j} (f_{k-1}) \right\|_k < 2^{-k}.\label{delta}
\end{align*}
By applying Proposition \ref{EndGame} with $f=f_{k-1}$ and $p=p_k$ to each operator $T_{i}$, there exist an integer $n_k > n_{k-1}$ and a polynomial $q_k$ of degree greater than $k$ such that conditions (a),(c),(d), and (e) in the proof of Theorem \ref{HypThm} are satisfied, along with the following condition:
\begin{enumerate}
    \item[(b)] for each integer $i$ with $1\leq i\leq k$, $\|T_{i}^{n_k}((q_k+1)f_{k-1}) - p_k \|_k < \ds 2^{-k}$. 
\end{enumerate}
One can now see that the rest of the argument for this generalization proceeds as in the proof of Theorem \ref{HypThm}, with some straightforward modifications. We omit the details here.

\section*{Acknowledgements}
Funding: This work was supported by the St. Olaf College Collaborative Undergraduate Research and Inquiry program.

%\section*{References}
%\printbibliography

\end{document}